\documentclass{article}%
\usepackage{amssymb}
\usepackage{amsmath}
\usepackage{amsfonts}
\usepackage{graphicx}%
\providecommand{\U}[1]{\protect \rule{.1in}{.1in}}
\newtheorem{theorem}{Theorem}

\newtheorem{proposition}[theorem]{Proposition}
\newtheorem{remark}[theorem]{Remark}

\begin{document}

\title{\textbf{A relative-geometric treatment}\\ \textbf{of ruled surfaces}}
\author{\textbf{Georg Stamou, Stylianos Stamatakis, Ioannis Delivos\bigskip}\\ \emph{Aristotle University of Thessaloniki}\\ \emph{Department of Mathematics}\\ \emph{GR-54124 Thessaloniki, Greece}\\ \emph{e-mail: stamoug@math.auth.gr}}
\date{}
\maketitle

\begin{abstract}
\noindent We consider relative normalizations of ruled surfaces with 
non-vanishing Gaussian curvature $K$ in the Euclidean space $%
\mathbb{R}
^{3},$ which are characterized by the support functions $^{\left(
\alpha \right)  }q=\left \vert K\right \vert ^{\alpha}$ for $\alpha \in%
\mathbb{R}
$. All ruled surfaces for which the relative normals, the Pick invariant or
the Tchebychev vector field have some specific properties are determined. We
conclude the paper by the study of the affine normal image of a non-conoidal
ruled surface.\smallskip

\smallskip


\end{abstract}

\section{Introduction}

Relatively normalized hypersurfaces with non-vanishing Gaussian curvature $K$
in the Euclidean space $%
\mathbb{R}
^{n+1}$, whose relative normalizations are characterized by the
\textit{support functions }$^{\left(  \alpha \right)  }q=\left \vert
K\right \vert ^{\alpha}$ for $\alpha \in%
\mathbb{R}
$ (see Section 2), have been studied in the last two decades by a number of
authors and many results have been derived. The one-parameter family of
relative normalizations $^{\left(  \alpha \right)  }\bar{y},$ which is
determined by the support functions $^{\left(  \alpha \right)  }q$ and was
introduced by F. Manhart \cite{Manhart5}, deserves special interest, because
in this family are contained among other relative normalizations the
\textit{Euclidean normalization }(for $\alpha=0$) as well as the
\textit{equiaffine normalization} (for $\alpha=1/\left(  n+2\right)  $). A
class of surfaces in $%
\mathbb{R}
^{3}$, relatively normalized by $^{\left(  \alpha \right)  }\bar{y}$, which has
been investigated in the past, is the one of ruled surfaces (see e.g.
\cite{Manhart3}, \cite{Manhart4}, \cite{Stamou}). These surfaces are further
discussed in the first part of the present work. The affine normal image of a
non-conoidal ruled surface is studied in the second part.

\section{Relative normalizations of surfaces}

In the Euclidean space $%
\mathbb{R}
^{3}$ let $\Phi:\bar{x}=\bar{x}(u^{1},u^{2}):U\subset%
\mathbb{R}
^{2}\longrightarrow%
\mathbb{R}
^{3}$ be an injective $C^{r}$-immersion with Gaussian curvature $K\neq0$
$\forall$ $(u^{1},u^{2})\in U.$ A $C^{s}$-mapping $\bar{y}:U\longrightarrow%
\mathbb{R}
^{3}$ ($r>s\geq1$) is called a $C^{s}$-\textit{relative normalization}
\footnote{For notations and definitions the reader is referred to
\cite{Schirokow} and \cite{Schneider}.} if
\begin{equation}
\bar{y}(\alpha)\notin T_{P}\Phi,\quad \frac{\partial \bar{y}}{\partial u^{i}%
}(\alpha)\in T_{P}\Phi \quad(i=1,2),\quad P=\bar{x}(\alpha), \tag{2.1}%
\label{21}%
\end{equation}
at every point $P\in \Phi$, where $T_{P}\Phi$ is the tangent vector space of
$\Phi$ in $P.\smallskip$

\noindent The \textit{covector }$\bar{X}$ of the tangent plane is defined by
\begin{equation}
\langle \bar{X},\frac{\partial \bar{x}}{\partial u^{i}}\rangle=0\quad
(i=1,2)\quad \text{and}\quad \langle \bar{X},\bar{y}\rangle=1, \tag{2.2}%
\label{22}%
\end{equation}
where $\langle~,~\rangle$\ denotes the standard scalar product in $%
\mathbb{R}
^{3}$. Using $\bar{X},$ the \textit{relative metric }$G$ is introduced by%
\begin{equation}
G_{ij}=\langle \bar{X},\frac{\partial^{2}\bar{x}}{\partial u^{i}\partial u^{j}%
}\rangle. \tag{2.3}\label{23}%
\end{equation}
From now on, we will use the tensor $G_{ij}$ for raising and lowering the
indices in the sense of classical tensor notation. Let $\bar{\xi
}:U\longrightarrow%
\mathbb{R}
^{3}$ be the \textit{Euclidean normalization }of $\Phi.$ The \textit{support
function }of the relative normalization $\bar{y}$ is defined by%
\begin{equation}
q:=\langle \bar{\xi},\bar{y}\rangle:U\longrightarrow%
\mathbb{R}
,\quad q\in C^{s}(U). \tag{2.4}\label{24}%
\end{equation}
By virtue of (\ref{21}), the support function $q$ is not vanishing in $U;$
moreover because of (\ref{22}) it is
\begin{equation}
\bar{X}=q^{-1}\bar{\xi}. \tag{2.5}\label{25}%
\end{equation}
From (\ref{23}) and (\ref{25}) we obtain%
\begin{equation}
G_{ij}=q^{-1}h_{ij}, \tag{2.6}\label{26}%
\end{equation}
where $h_{ij}$ are the components of the second fundamental form of $\Phi$. We
mention that given a support function $q$ the relative normalization $\bar{y}$
is uniquely determined and possesses the following parametrization (see
\cite[p. 197]{Manhart5})%
\begin{equation}
\bar{y}=-h^{(ij)}\frac{\partial q}{\partial u^{i}}\frac{\partial \bar{x}%
}{\partial u^{j}}+q\bar{\xi}, \tag{2.7}\label{27}%
\end{equation}
where $h^{(ij)}$ are the components of the inverse tensor of $h_{ij}%
$.\smallskip

\noindent Let $^{G}\nabla_{i}f$ denote the covariant derivative with respect
to $G$ of a differentiable function $f,$ which is defined on $\Phi.$ The
symmetric \textit{Darboux tensor }is defined by%
\begin{equation}
A_{jkl}:=\langle \bar{X},^{G}\nabla_{l}^{G}\nabla_{k}\frac{\partial \bar{x}%
}{\partial u^{j}}\rangle, \tag{2.8}\label{28}%
\end{equation}
the \textit{Tchebychev vector }by%
\begin{equation}
T_{i}=\frac{1}{2}A_{ijk}G^{jk}=\frac{1}{2}A_{ij}^{j}, \tag{2.9}\label{29}%
\end{equation}
and the \textit{Pick invariant }by%
\begin{equation}
J:=\frac{1}{2}A_{jkl}A^{jkl}. \tag{2.10}\label{210}%
\end{equation}

\section{Relatively normalized ruled surfaces by $^{(\alpha)}\bar{y}$}

\textbf{3.1. }Let $\Phi \subset%
\mathbb{R}
^{3}$ be a skew ruled $C^{2}$-surface. We denote by $\bar{s}(u),$ $u\in I$
($I\subset%
\mathbb{R}
$ open interval) the position vector of the line of striction of $\Phi$\ and
by $\bar{e}(u)$ the unit vector pointing along the rulings. Moreover we can
choose the parameter $u$ to be the arc length along the spherical curve
$\bar{e}(u)$. Then a parametrization of the ruled surface $\Phi$ over the
region $U:=I\times%
\mathbb{R}
$ of the $\left(  u,v\right)  $-plane is%
\begin{equation}
\bar{x}(u,v)=\bar{s}(u)+v\bar{e}(u)\quad \text{with}\quad \left \vert \bar
{e}\right \vert =|\bar{e}\,%
\acute{}%
\,|=1,\text{\quad}\langle \bar{s}\,%
\acute{}%
\,(u),\bar{e}\,%
\acute{}%
\,(u)\rangle=0\text{ in }I, \tag{3.1}\label{301}%
\end{equation}
where the prime denotes differentiation with respect to $u$. The moving frame
of $\Phi,$ consisting of the vector $\bar{e}(u)$, the central normal vector
$\bar{n}:=\bar{e}$\thinspace$%
\acute{}%
$ and the central tangent vector $\bar{z}:=\bar{e}\times \bar{n},$ moves along
the line of striction and satisfies the relations \cite[p. 62f]{Kruppa}%
\begin{equation}
\bar{e}\,%
\acute{}%
=\bar{n},\quad \bar{n}\,%
\acute{}%
=-\bar{e}+\kappa \bar{z},\quad \bar{z}\,%
\acute{}%
=-\kappa \bar{n}, \tag{3.2}\label{302}%
\end{equation}
where $\kappa=(\bar{e},\bar{e}$\thinspace$%
\acute{}%
,\bar{e}$\thinspace$%
\acute{}%
$\thinspace$%
\acute{}%
$\thinspace$)$ denotes the \textit{conical curvature }of $\Phi.$ Consider the
\textit{parameter of distribution }$\delta=(\bar{s}$\thinspace$%
\acute{}%
,\bar{e},\bar{e}$\thinspace$%
\acute{}%
$\thinspace$)$ and the \textit{striction}\emph{ }$\sigma:=\sphericalangle
(\bar{e},\bar{s}$\thinspace$%
\acute{}%
$\thinspace$)$ ($-\frac{\pi}{2}<\sigma \leq \frac{\pi}{2}$,
$\operatorname*{sign}\sigma=\operatorname*{sign}\delta).$ Then the tangent
vector $\bar{s}\,%
\acute{}%
$ of the line of striction has the expression%
\begin{equation}
\bar{s}\,%
\acute{}%
=\delta \left(  \lambda \bar{e}+\bar{z}\right)  \quad \text{with\quad}%
\lambda:=\cot \sigma. \tag{3.3}\label{303}%
\end{equation}
When the invariants $\kappa(u),\delta(u)$ and $\lambda(u)$ (fundamental
invariants) are given, then there exists up to rigid motions of the space $%
\mathbb{R}
^{3}$ a unique ruled surface $\Phi$, whose fundamental invariants are the
given. The components $g_{ij}$ and $h_{ij}$ of the first and the second
fundamental tensors in the coordinates $u^{1}:=u$, $u^{2}:=v$ are the
following%
\begin{equation}
{\normalsize (g}_{ij}{\normalsize )=}\left(
\begin{array}
[c]{cc}%
v^{2}+\delta^{2}\left(  \lambda^{2}+1\right)  & \delta \lambda \\
\delta \lambda & 1
\end{array}
\right)  , \tag{3.4}\label{304}%
\end{equation}%
\begin{equation}
{\small (h}_{ij}{\small )=}\frac{1}{w}\left(
\begin{array}
[c]{cc}%
-\left[  \kappa v^{2}{\small +}\delta \,%
\acute{}%
\,v{\small +}\delta^{2}\left(  \kappa-\lambda \right)  \right]  & \delta \\
\delta & 0
\end{array}
\right)  , \tag{3.5}\label{305}%
\end{equation}
where $w:=\sqrt{v^{2}+\delta^{2}}$. For the Gaussian curvature $K$ of $\Phi$
the following relation holds%
\begin{equation}
K=-\frac{\delta^{2}}{w^{4}}. \tag{3.6}\label{306}%
\end{equation}

We consider now the relative normalizations $^{\left(  \alpha \right)  }\bar
{y}:U\longrightarrow%
\mathbb{R}
^{3}$, which are introduced by F. Manhart \cite{Manhart5} and, on account of
(\ref{27}), are defined by the support functions%
\begin{equation}
^{\left(  \alpha \right)  }q:=|K|^{\alpha},\quad \alpha \in%
\mathbb{R}
. \tag{3.7}\label{307}%
\end{equation}
We denote by $^{\left(  \alpha \right)  }J$ the Pick invariant. In order to
compute it we firstly have from (\ref{28}) the relation \cite[p.
196]{Manhart5}%
\begin{equation}
A_{jkl}=\frac{1}{q}\langle \bar{\xi},\frac{\partial^{3}\bar{x}}{\partial
u^{j}\partial u^{k}\partial u^{l}}\rangle-\frac{1}{2}\left(  \frac{\partial
G_{jk}}{\partial u^{l}}+\frac{\partial G_{kl}}{\partial u^{j}}+\frac{\partial
G_{lj}}{\partial u^{k}}\right)  . \tag{3.8}\label{308}%
\end{equation}
Using (\ref{308}), we find for the ruled surface (\ref{301}), which is
relatively normalized by $^{\left(  \alpha \right)  }\bar{y},$%
\begin{equation}
A_{112}=\frac{\left(  4\alpha-1\right)  |\delta|^{-2\alpha}}{2w^{3-4\alpha}%
}\left[  \kappa v^{3}+2\delta \,%
\acute{}%
\,v^{2}+\delta^{2}\left(  \kappa-\lambda \right)  v-\delta^{2}\delta \,%
\acute{}%
\, \right]  , \tag{3.9}\label{309}%
\end{equation}%
\begin{equation}
A_{221}=\frac{\varepsilon \left(  1-4\alpha \right)  |\delta|^{-2\alpha+1}%
}{w^{3-4\alpha}}v,\quad A_{222}=0, \tag{3.10}\label{310}%
\end{equation}%
\begin{align}
A_{111}  &  =\frac{\varepsilon|\delta|^{-2\alpha-1}}{2w^{3-4\alpha}}\{ \left(
\delta \kappa \,%
\acute{}%
-6\alpha \delta \,%
\acute{}%
\, \kappa \right)  v^{4}\tag{3.11}\label{311}\\
&  +\left[  -2\delta^{2}\left(  1+\kappa \lambda \right)  +\delta \delta \,%
\acute{}%
\,%
\acute{}%
-6\alpha \delta \,%
\acute{}%
\,^{2}\right]  v^{3}\nonumber \\
&  +\delta^{2}\left[  \delta \left(  2\kappa \,%
\acute{}%
+\lambda \,%
\acute{}%
\, \right)  -\delta \,%
\acute{}%
\, \kappa+2\left(  3\alpha-1\right)  \delta \,%
\acute{}%
\, \lambda \right]  v^{2}\nonumber \\
&  +\delta^{2}\left[  -2\delta^{2}\left(  1+\kappa \lambda \right)
+\delta \delta \,%
\acute{}%
\,%
\acute{}%
+3\left(  2\alpha-1\right)  \delta \,%
\acute{}%
\,^{2}\right]  v\nonumber \\
&  +\delta^{4}\left[  \left(  6\alpha-1\right)  \left(  \kappa-\lambda \right)
\delta \,%
\acute{}%
+\delta \left(  \kappa \,%
\acute{}%
+\lambda \,%
\acute{}%
\, \right)  \right]  \},\nonumber
\end{align}
where $\varepsilon=\operatorname*{sign}\delta.$ Furthermore, since the tensor
$A_{jkl}$ is symmetric, we have%
\begin{equation}
A_{112}=A_{211}=A_{121},\quad A_{221}=A_{212}=A_{122}. \tag{3.12}\label{312}%
\end{equation}
The components $A^{jkl}$ can be found from $A^{jkl}=G^{ji}G^{km}G^{lr}A_{imr}%
$. Inserting $A_{jkl}$ and $A^{jkl}$ in (\ref{210}), it turns out that%
\begin{equation}
^{\left(  \alpha \right)  }J=\frac{3\left(  4\alpha-1\right)  ^{2}%
|\delta|^{2\alpha-2}}{2w^{4\alpha+3}}\left[  \kappa v^{4}+\delta^{2}\left(
\kappa-\lambda \right)  v^{2}+\delta^{2}\delta \,%
\acute{}%
\,v\right]  . \tag{3.13}\label{313}%
\end{equation}
Let $^{\left(  \alpha \right)  }\bar{T}$ be the corresponding Tchebychev vector
of $^{\left(  \alpha \right)  }\bar{y}$. Applying similar computations as above
we can find $^{\left(  \alpha \right)  }\bar{T}$ as follows: In view of
$T^{i}=G^{ij}T_{j}$ and using (\ref{29}) we obtain%
\begin{equation}
T^{1}=\frac{\varepsilon \left(  1-4\alpha \right)  |\delta|^{2\alpha-1}%
}{w^{4\alpha+1}}v, \tag{3.14}\label{314}%
\end{equation}%
\begin{equation}
T^{2}=\frac{\left(  1-4\alpha \right)  |\delta|^{2\alpha-2}}{2w^{4\alpha+1}%
}\left[  2\kappa v^{3}+\delta \,%
\acute{}%
\,v^{2}+2\delta^{2}\left(  \kappa-\lambda \right)  v+\delta^{2}\delta \,%
\acute{}%
\, \right]  . \tag{3.15}\label{315}%
\end{equation}
Then, substituting (\ref{314}) and (\ref{315}) into%
\begin{equation}
^{\left(  \alpha \right)  }\bar{T}=T^{1}\frac{\partial \bar{x}}{\partial
u}+T^{2}\frac{\partial \bar{x}}{\partial v}, \tag{3.16}\label{316}%
\end{equation}
we obtain%
\begin{equation}
^{\left(  \alpha \right)  }\bar{T}\!=\! \frac{\left(  1\!-\!4\alpha \right)
|\delta|^{2\alpha-2}}{2w^{4\alpha+1}}\left[  \left(  2\kappa v^{3}\!+\!
\delta \,%
\acute{}%
\,v^{2}\!+\!2\delta^{2}\kappa v+\delta^{2}\delta \,%
\acute{}%
\, \right)  \bar{e}\!+\!2\delta v^{2}\bar{n}\!+\!2\delta^{2}v\bar{z}\right]  .
\tag{3.17}\label{317}%
\end{equation}
In the following paragraphs we will discuss questions on ruled surfaces
relatively normalized by $^{\left(  \alpha \right)  }\bar{y}$, which are
related with the relative normals, the Pick invariant and the Tchebychev
vector field.\bigskip

\noindent \textbf{3.2.} In this paragraph we treat ruled surfaces, whose
relative normal vectors $^{\left(  \alpha \right)  }\bar{y}$ are parallel to a
fixed plane $E$. The vectors $^{\left(  \alpha \right)  }\bar{y}$ are given by
the relation (see \cite[p. 212]{Stamou})%
\begin{equation}
^{\left(  \alpha \right)  }\bar{y}=A_{1}\bar{e}+A_{2}\bar{n}+A_{3}\bar{z},
\tag{3.18}\label{318}%
\end{equation}
where%
\begin{equation}
A_{1}=\frac{2\alpha \left(  2\kappa v+\delta \,%
\acute{}%
\, \right)  |\delta|^{2\alpha-2}}{w^{4\alpha-1}}, \tag{3.19}\label{319}%
\end{equation}%
\begin{equation}
A_{2}=\frac{\varepsilon \left(  4\alpha v^{2}+\delta^{2}\right)  |\delta
|^{2\alpha-1}}{w^{4\alpha+1}}, \tag{3.20}\label{320}%
\end{equation}%
\begin{equation}
A_{3}=\frac{\left(  4\alpha-1\right)  |\delta|^{2\alpha}}{w^{4\alpha+1}}v.
\tag{3.21}\label{321}%
\end{equation}
Let $\bar{c}\neq \bar{0}$ be a constant normal vector of the plane $E$. Because
of (\ref{319})-(\ref{321}) the assumption $\langle^{\left(  \alpha \right)
}\bar{y},\bar{c}\rangle=0$ leads to the relation%
\begin{align}
&  4\alpha \kappa \langle \bar{e},\bar{c}\rangle v^{3}+2\alpha \left(  \delta \,%
\acute{}%
\, \langle \bar{e},\bar{c}\rangle+2\delta \langle \bar{n},\bar{c}\rangle \right)
v^{2}\tag{3.22}\label{322}\\
&  +\delta^{2}\left[  4\alpha \kappa \langle \bar{e},\bar{c}\rangle+\left(
4\alpha-1\right)  \langle \bar{z},\bar{c}\rangle \right]  v+2\alpha \delta
^{2}\delta \,%
\acute{}%
\, \langle \bar{e},\bar{c}\rangle+\delta^{3}\langle \bar{n},\bar{c}%
\rangle=0.\nonumber
\end{align}
The polynomial on the left hand side of (\ref{322}) is of degree three in $v$
and vanishes for all $u\in I$ and infinite values for $v\in%
\mathbb{R}
$. Comparing its coefficients with those of the zero polynomial we obtain%
\begin{equation}
\alpha \kappa \langle \bar{e},\bar{c}\rangle=0, \tag{3.23}\label{323}%
\end{equation}%
\begin{equation}
\alpha \left(  \delta \,%
\acute{}%
\, \langle \bar{e},\bar{c}\rangle+2\delta \langle \bar{n},\bar{c}\rangle \right)
=0, \tag{3.24}\label{324}%
\end{equation}%
\begin{equation}
4\alpha \kappa \langle \bar{e},\bar{c}\rangle+\left(  4\alpha-1\right)
\langle \bar{z},\bar{c}\rangle=0, \tag{3.25}\label{325}%
\end{equation}%
\begin{equation}
2\alpha \delta \,%
\acute{}%
\, \langle \bar{e},\bar{c}\rangle+\delta \langle \bar{n},\bar{c}\rangle=0.
\tag{3.26}\label{326}%
\end{equation}
If $\alpha \neq1/4$, then from (\ref{323})-(\ref{326}) we get $\langle \bar
{n},\bar{c}\rangle=\langle \bar{z},\bar{c}\rangle=0,$ i.e. $\bar{e}/\!/\bar{c}%
$, which is impossible. If $\alpha=1/4,$ then according to \cite{Opozda}
$\Phi$ is a conoidal surface ($\kappa=0$). The same result follows also from
(\ref{323})-(\ref{326}). The above discussion gives rise to the following

\begin{proposition}
Let $\Phi$ be a ruled $C^{3}$-surface $\Phi,$ free of torsal rulings, which is
relatively normalized by $^{\left(  \alpha \right)  }\bar{y}$ ($\alpha \in%
\mathbb{R}
$). If the relative normals of $\Phi$ are parallel to a fixed plane, then
$\alpha=1/4$ and $\Phi$ is a conoidal surface.\medskip
\end{proposition}

\noindent \textbf{3.3.} In this paragraph we classify the ruled surfaces
$\Phi \subset%
\mathbb{R}
^{3}$, whose Tchebychev vectors $^{\left(  \alpha \right)  }\bar{T}$
$(\alpha \neq1/4)$ are tangent to some geometrically distinguished families of
curves of $\Phi.\smallskip$

\noindent A first result in this direction is obtained immediately from
(\ref{317}): \textit{The vectors }$^{\left(  \alpha \right)  }\bar{T}$
\textit{are tangent to the orthogonal trajectories of the rulings, if and only
if }$\kappa=\delta$\thinspace$%
\acute{}%
=0,$ \textit{i.e. if and only if }$\Phi$ \textit{is a conoidal surface with
constant parameter of distribution.\smallskip}

\noindent We consider a directrix $\Gamma$ of $\Phi$ defined by $v=v(u)$. In
view of (\ref{317}) we find that the vectors $^{\left(  \alpha \right)  }%
\bar{T}$ along $\Gamma$ are \smallskip

$\bullet \quad$tangent to $\Gamma$ if and only if
\begin{equation}
2\kappa v^{3}+\delta \,%
\acute{}%
\,v^{2}+2\delta \left[  \delta \left(  \kappa-\lambda \right)  -v\,%
\acute{}%
\, \right]  v+\delta^{2}\delta \,%
\acute{}%
=0,\text{ and} \tag{3.27}\label{327}%
\end{equation}

$\bullet \quad$orthogonal to $\Gamma$ if and only if%
\begin{equation}
\left(  \delta \lambda+v\,%
\acute{}%
\, \right)  \left(  2\kappa v+\delta \,%
\acute{}%
\, \right)  +2\delta v=0. \tag{3.28}\label{328}%
\end{equation}
Furthermore we consider the following curves of $\Phi$:\textit{\smallskip}

a) the curved asymptotic lines,

b) the curves of constant striction distance $(u$-curves) and

c) the $K$-curves, i.e. the curves along which the Gaussian curvature is con-

stant \cite{Sachs}.\textit{\smallskip}

\noindent The corresponding differential equations of these families of curves
are%
\begin{equation}
\kappa v^{2}+\delta \,%
\acute{}%
\,v+\delta^{2}\left(  \kappa-\lambda \right)  -2\delta v\,%
\acute{}%
=0, \tag{3.29}\label{329}%
\end{equation}%
\begin{equation}
v\,%
\acute{}%
=0, \tag{3.30}\label{330}%
\end{equation}%
\begin{equation}
2\delta vv\,%
\acute{}%
+\delta \,%
\acute{}%
\left(  \delta^{2}-v^{2}\right)  =0, \tag{3.31}\label{331}%
\end{equation}
respectively. On account of (\ref{327}), we find that the vectors $^{\left(
\alpha \right)  }\bar{T}$ ($\alpha \neq1/4)$ are tangent to one of these
families of curves if and only if
\begin{equation}
\kappa v^{3}+\delta^{2}\left(  \kappa-\lambda \right)  v+\delta^{2}\delta \,%
\acute{}%
=0, \tag{3.32}\label{332}%
\end{equation}%
\begin{equation}
2\kappa v^{3}+\delta \,%
\acute{}%
\,v^{2}+2\delta^{2}\left(  \kappa-\lambda \right)  v+\delta^{2}\delta \,%
\acute{}%
=0, \tag{3.33}\label{333}%
\end{equation}%
\begin{equation}
\kappa v^{3}+\delta^{2}(\kappa-\lambda)v+\delta^{2}\delta \,%
\acute{}%
\,=0, \tag{3.34}\label{334}%
\end{equation}
respectively. Since each condition is satisfied for all $(u,v)\in U$, we
obtain $\kappa=\lambda=\delta$\thinspace$%
\acute{}%
=0.$ Therefore $\Phi$ lies on a right helicoid. The same result arises, when
we suppose that the Pick invariant $^{\left(  \alpha \right)  }J$ ($\alpha
\neq1/4)$ vanishes identically, as it follows immediately by means of
(\ref{313}). The above constitute the proof of

\begin{proposition}
Let $\Phi \subset%
\mathbb{R}
^{3}$ be a ruled $C^{3}$-surface, free of torsal rulings, which is relatively
normalized by $^{\left(  \alpha \right)  }\bar{y}$ ($\alpha \in%
\mathbb{R}
/\left \{  1/4\right \}  $). The following properties are equivalent:\newline(i)
The Tchebychev vectors $^{\left(  \alpha \right)  }\bar{T}$ are tangent to the
curved asymptotic lines or to the curves of constant striction distance or to
the $K$-curves.\newline(ii) The Pick invariant vanishes identically.\newline%
(iii) $\Phi$ lies on a right helicoid.
\end{proposition}

\noindent Continuing this line of work, we require now that the vectors
$^{\left(  \alpha \right)  }\bar{T}$ are tangent to the orthogonal trajectories
of the $u$-curves or to the $K$-curves. On account of (\ref{330}) and
(\ref{331}) and by virtue of (\ref{328}), we obtain in the first case the
condition%
\begin{equation}
2\left(  1+\kappa \lambda \right)  v+\delta \,%
\acute{}%
\, \lambda=0 \tag{3.35}\label{335}%
\end{equation}
and in the second case the condition%
\begin{equation}
2\delta \,%
\acute{}%
\, \kappa v^{3}+\left[  \delta \,%
\acute{}%
\,^{2}+4\delta^{2}\left(  1+\kappa \lambda \right)  \right]  v^{2}-2\delta
^{2}\delta \,%
\acute{}%
\left(  \kappa-\lambda \right)  v-\delta^{2}\delta \,%
\acute{}%
\,^{2}=0, \tag{3.36}\label{336}%
\end{equation}
which are satisfied for every $(u,v)\in U$. Thus, we have%
\begin{equation}
\delta \,%
\acute{}%
=1+\kappa \lambda=0. \tag{3.37}\label{337}%
\end{equation}
Next, we assume that the vectors $^{\left(  \alpha \right)  }\bar{T}$ are
tangent to one family of (Euclidean) lines of curvature. We substitute the
derivative $v$\thinspace$%
\acute{}%
$\thinspace$(u)$ from (\ref{327}) into the differential equation of the lines
of curvature%
\begin{equation}
g_{12}h_{11}-g_{11}h_{12}+\left(  g_{22}h_{11}-g_{11}h_{22}\right)  v\,%
\acute{}%
+\left(  g_{22}h_{12}-g_{12}h_{22}\right)  v\,%
\acute{}%
\,^{2}=0 \tag{3.38}\label{338}%
\end{equation}
and taking into account (\ref{304}) and (\ref{305}), we obtain once more
(\ref{336}). Hence we get again the conditions (\ref{337}), which show that
$\Phi$ is a ruled surface with constant parameter of distribution and whose
line of striction is a line of curvature. This characterizes the
\textit{Edlinger surfaces }(see \cite[p. 31]{Hoschek}) which, by definition,
are ruled surfaces whose osculating quadrics are rotational hyperboloids (see
\cite[p. 36]{Hoschek}). Thus the following proposition is valid:

\begin{proposition}
Let $\Phi \subset%
\mathbb{R}
^{3}$ be a ruled $C^{3}$-surface, free of torsal rulings, which is relatively
normalized by $^{\left(  \alpha \right)  }\bar{y}$ $(\alpha \in%
\mathbb{R}
/\left \{  1/4\right \}  ).$ Assume that the Tchebychev vectors $^{\left(
\alpha \right)  }\bar{T}$ are tangent to the orthogonal trajectories of the
curves of constant striction distance or to the $K$-curves or tangent to one
family of (Euclidean) lines of curvature. Then $\Phi$ is an Edlinger surface.

\begin{remark}
(a) On a right helicoid the three families of curves, which are mentioned in
Proposition 2(i), coincide.\newline(b) On an Edlinger surface the curves of
constant striction distance and the $K$-curves coincide and they are lines of
curvature.\medskip
\end{remark}
\end{proposition}

\noindent \textbf{3.4.} In this paragraph we will study ruled surfaces, whose
divergence or rotation (curl), with respect to the metric $g_{ij}du^{i}%
du^{j},$ of the vector field $^{\left(  \alpha \right)  }\bar{T}$ ($\alpha
\neq1/4)$ vanishes identically.\textit{\smallskip}

The divergence of $^{\left(  \alpha \right)  }\bar{T}(u,v)$ is given by the
relation \cite[p. 121]{Strubecker}%
\begin{equation}
\operatorname{div}\left(  ^{\left(  \alpha \right)  }\bar{T}\left(  u,v\right)
\right)  =\frac{1}{w}\left(  \frac{\partial \left(  wT^{1}\right)  }{\partial
u}+\frac{\partial \left(  wT^{2}\right)  }{\partial v}\right)  . \tag{3.39}%
\label{340}%
\end{equation}
After a short calculation we obtain%
\begin{align}
&  \operatorname{div}\left(  ^{\left(  \alpha \right)  }\bar{T}\left(
u,v\right)  \right)  =\frac{\left(  1-4\alpha \right)  \left \vert
\delta \right \vert ^{2\alpha-2}}{w^{4\alpha+3}}\{ \left(  3-4\alpha \right)
\kappa v^{4}\tag{3.40}\label{341}\\
&  +\delta^{2}\left[  4\left(  1-\alpha \right)  \kappa+\left(  4\alpha
-1\right)  \lambda \right]  v^{2}-4\alpha \delta^{2}\delta \,%
\acute{}%
\,v+\delta^{4}\left(  \kappa-\lambda \right)  \}.\nonumber
\end{align}
The identical vanishing of $\operatorname{div}\left(  ^{\left(  \alpha \right)
}\bar{T}\left(  u,v\right)  \right)  $ implies the following conditions:%
\begin{equation}
\left(  3-4\alpha \right)  \kappa=0, \tag{3.41}\label{342}%
\end{equation}%
\begin{equation}
4\left(  1-\alpha \right)  \kappa+\left(  4\alpha-1\right)  \lambda=0,
\tag{3.42}\label{343}%
\end{equation}%
\begin{equation}
\alpha \delta \,%
\acute{}%
=0, \tag{3.43}\label{344}%
\end{equation}%
\begin{equation}
\kappa-\lambda=0. \tag{3.44}\label{345}%
\end{equation}
Elementary treatment of the above system yields: a) if $\alpha=0,~$then
$\kappa=\lambda=0$, i.e. $\Phi$ is a right conoid, b) if $\alpha \in%
\mathbb{R}
/\left \{  0,1/4\right \}  $, then $\kappa=\lambda=\delta$\thinspace$%
\acute{}%
=0,$ which means that $\Phi$ lies on a right helicoid.\textit{\smallskip}

We compute now the rotation of $^{\left(  \alpha \right)  }\bar{T}(u,v)$.
According to \cite[p. 125]{Strubecker} holds%
\begin{equation}
\operatorname*{rot}\left(  ^{\left(  \alpha \right)  }\bar{T}\left(
u,v\right)  \right)  =\frac{1}{w}\left[  \frac{\partial \left(  T^{1}%
g_{12}+T^{2}g_{22}\right)  }{\partial u}-\frac{\partial \left(  T^{1}%
g_{11}+T^{2}g_{12}\right)  }{\partial v}\right]  . \tag{3.45}\label{346}%
\end{equation}
On account of (\ref{304}), (\ref{314}) and (\ref{315}) we then find
\begin{align}
\operatorname*{rot}\left(  ^{\left(  \alpha \right)  }\bar{T}\left(
u,v\right)  \right)   &  =\frac{\varepsilon \left(  1-4\alpha \right)
\left \vert \delta \right \vert ^{2\alpha-3}}{2w^{4\alpha+4}}\{ \left[  4\left(
\alpha-1\right)  \delta \,%
\acute{}%
\, \kappa+2\delta \kappa \,%
\acute{}%
\, \right]  v^{5}\tag{3.46}\label{347}\\
&  +\left[  2\left(  \alpha-1\right)  \delta \,%
\acute{}%
\,^{2}+\delta \delta \,%
\acute{}%
\,%
\acute{}%
+4\delta^{2}\left(  2\alpha-1\right)  \left(  1+\kappa \lambda \right)  \right]
v^{4}\nonumber \\
&  +\delta^{2}\left[  4\delta \kappa \,%
\acute{}%
-6\delta \,%
\acute{}%
\, \kappa+\left(  4\alpha-1\right)  \delta \,%
\acute{}%
\, \lambda \right]  v^{3}\nonumber \\
&  +\delta^{2}\left[  -3\delta \,%
\acute{}%
\,^{2}+2\delta \delta \,%
\acute{}%
\,%
\acute{}%
+2\delta^{2}\left(  4\alpha-3\right)  \left(  1+\kappa \lambda \right)  \right]
v^{2}\nonumber \\
&  +\delta^{4}\left[  -2\left(  2\alpha+1\right)  \delta \,%
\acute{}%
\, \kappa+2\delta \kappa \,%
\acute{}%
+\left(  4\alpha-1\right)  \delta \,%
\acute{}%
\, \lambda \right]  v\nonumber \\
&  +\delta^{4}\left[  -\left(  2\alpha+1\right)  \delta \,%
\acute{}%
\,^{2}+\delta \delta \,%
\acute{}%
\,%
\acute{}%
-2\delta^{2}\left(  1+\kappa \lambda \right)  \right]  \},\nonumber
\end{align}
which vanishes identically if and only if%
\begin{equation}
2\left(  \alpha-1\right)  \delta \,%
\acute{}%
\, \kappa+\delta \kappa \,%
\acute{}%
=0, \tag{3.47}\label{348}%
\end{equation}%
\begin{equation}
2\left(  \alpha-1\right)  \delta \,%
\acute{}%
\,^{2}+\delta \delta \,%
\acute{}%
\,%
\acute{}%
+4\delta^{2}\left(  2\alpha-1\right)  \left(  1+\kappa \lambda \right)  =0,
\tag{3.48}\label{349}%
\end{equation}%
\begin{equation}
4\delta \kappa \,%
\acute{}%
-6\delta \,%
\acute{}%
\, \kappa+\left(  4\alpha-1\right)  \delta \,%
\acute{}%
\, \lambda=0, \tag{3.49}\label{350}%
\end{equation}%
\begin{equation}
-3\delta \,%
\acute{}%
\,^{2}+2\delta \delta \,%
\acute{}%
\,%
\acute{}%
+2\delta^{2}\left(  4\alpha-3\right)  \left(  1+\kappa \lambda \right)  =0,
\tag{3.50}\label{351}%
\end{equation}%
\begin{equation}
-2\left(  2\alpha+1\right)  \delta \,%
\acute{}%
\, \kappa+2\delta \kappa \,%
\acute{}%
+\left(  4\alpha-1\right)  \delta \,%
\acute{}%
\, \lambda=0, \tag{3.51}\label{352}%
\end{equation}%
\begin{equation}
-\left(  2\alpha+1\right)  \delta \,%
\acute{}%
\,^{2}+\delta \delta \,%
\acute{}%
\,%
\acute{}%
-2\delta^{2}\left(  1+\kappa \lambda \right)  =0. \tag{3.52}\label{353}%
\end{equation}
From (\ref{348}) and (\ref{350}) we obtain
\begin{equation}
\delta \,%
\acute{}%
\left(  2\kappa-\lambda \right)  =0 \tag{3.53}\label{354}%
\end{equation}
and from (\ref{351}) and (\ref{353})%
\begin{equation}
\delta \delta \,%
\acute{}%
\,%
\acute{}%
-2\alpha \delta \,%
\acute{}%
\,^{2}=0. \tag{3.54}\label{355}%
\end{equation}
From (\ref{353}) and (\ref{355}), and if $\lambda=2\kappa,$ we deduce that%
\[
\delta \,%
\acute{}%
\,^{2}+2\delta^{2}+4\delta^{2}\kappa^{2}=0,
\]
a contradiction, since $\delta \neq0$. From (\ref{348}) and for $\delta
$\thinspace$%
\acute{}%
=0$ we then obtain $\kappa$\thinspace$%
\acute{}%
=0$ and from (\ref{353}) $1+\kappa \lambda=0.$ Thus, $\Phi$ is an Edlinger
surface, whose osculating hyperboloids are congruent \cite[p. 36]{Hoschek}.
The above discussion can be summarized in

\begin{proposition}
Let $\Phi \subset%
\mathbb{R}
^{3}$ be a ruled $C^{3}$-surface, free of torsal rulings, which is relatively
normalized by $^{\left(  \alpha \right)  }\bar{y}$ $(\alpha \in%
\mathbb{R}
/\left \{  1/4\right \}  ).$ For the Tchebychev vector field $^{\left(
\alpha \right)  }\bar{T}$ the following properties are valid:\newline(i)
$\operatorname{div}^{\left(  0\right)  }\bar{T}\equiv0,$ when $\Phi$ is a
right conoid.\newline(ii) $\operatorname{div}^{\left(  \alpha \right)  }\bar
{T}\equiv0$ $(\alpha \in%
\mathbb{R}
/\left \{  0,1/4\right \}  ),$ when $\Phi$ lies on a right helicoid.\newline%
(iii) $\operatorname*{rot}^{\left(  \alpha \right)  }\bar{T}\equiv0$
$(\alpha \in%
\mathbb{R}
/\left \{  1/4\right \}  ),$ when $\Phi$ is an Edlinger surface with congruent
osculating hyperboloids.
\end{proposition}

\section{The affine normal image of a ruled surface}

We consider a ruled surface $\Phi$ with the parametrization (\ref{301})$.$ A
parametrization of the affine normal image $\Phi^{\ast}$ of $\Phi$ is obtained
setting in (\ref{318}) $\alpha=1/4$:%
\begin{equation}
\bar{x}^{\ast}:=^{\left(  1/4\right)  }\bar{y}=\varepsilon \left \vert
\delta \right \vert ^{-1/2}\bar{n}+\frac{\left(  2\kappa v+\delta \,%
\acute{}%
\, \right)  }{2}\left \vert \delta \right \vert ^{-3/2}\bar{e},\quad
\varepsilon=\operatorname*{sign}\delta. \tag{4.1}\label{401}%
\end{equation}
From the above parametrization we have: \textit{Two ruled surfaces }%
$\Phi,\widetilde{\Phi},$ \textit{parametrized by (\ref{301}), with parallel
rulings and the same parameter of distribution have the same affine normal
image.\smallskip \ }

\noindent Hereafter, we consider only non-conoidal ($\kappa \neq0$) ruled
surfaces. In this case $\Phi^{\ast}$ is a ruled surface, whose rulings are
parallel to the corresponding rulings of $\Phi$. Its directrix%
\begin{equation}
\Gamma:\bar{r}^{\ast}=\varepsilon \left \vert \delta \right \vert ^{-1/2}\bar{n}
\tag{4.2}\label{402}%
\end{equation}
is a geometrically distinguished curve of $\Phi^{\ast}$: it is the sectional
curve of the director-cone of the surface of the central normals of $\Phi$
with $\Phi^{\ast}$. We can easily find the line of striction of $\Phi^{\ast}$
to be
\begin{equation}
\bar{s}^{\ast}=\frac{\delta \,%
\acute{}%
\, \left \vert \delta \right \vert ^{-3/2}}{2}\bar{e}+\varepsilon \left \vert
\delta \right \vert ^{-1/2}\bar{n}. \tag{4.3}\label{403}%
\end{equation}
We observe that setting in (\ref{401}) $v=0,$ we find the above vector. So we
have: \textit{The lines of striction of }$\Phi$ \textit{and} $\Phi^{\ast}$
\textit{correspond to each other}. Furthermore we conclude: \textit{The line
of striction of }$\Phi^{\ast}$ \textit{coincides with the curve }$\Gamma$
\textit{if and only if the parameter of distribution of }$\Phi$\textit{\ is
constant.\smallskip}

\noindent After a short computation we find as fundamental invariants
($\kappa^{\ast},\delta^{\ast},\lambda^{\ast}$) of $\Phi^{\ast}:$%
\begin{equation}
\kappa^{\ast}=\kappa, \tag{4.4}\label{404}%
\end{equation}%
\begin{equation}
\delta^{\ast}=\varepsilon \kappa \left \vert \delta \right \vert ^{-1/2},
\tag{4.5}\label{405}%
\end{equation}%
\begin{equation}
\lambda^{\ast}=\frac{2\delta \delta \,%
\acute{}%
\,%
\acute{}%
-3\delta \,%
\acute{}%
\,^{2}-4\delta^{2}}{4\delta^{2}\kappa}. \tag{4.6}\label{406}%
\end{equation}
From these relations we obtain:\textit{\smallskip}

\noindent a) $\Phi^{\ast}$ \textit{is an orthoid ruled surface }%
($\lambda^{\ast}=0$) \textit{if and only if }$\delta=1/\left(  c_{1}\sin
u+c_{2}\cos u\right)  ^{2}$, \textit{where} $c_{1},c_{2}=const.$%
\textit{\smallskip}

\noindent b) \textit{The line of striction of }$\Phi^{\ast}$\textit{ is a line
of curvature }($1+\kappa^{\ast}\lambda^{\ast}=0$)\textit{ if and only if
}$\delta=const.\neq0$\textit{ or }$\delta=c_{1}/u^{2},$ \textit{where}
$c_{1}=const.\neq0.$\textit{\smallskip}

\noindent The ruled surfaces $\Phi$ and $\Phi^{\ast}$ are congruent if and
only if $\kappa=\kappa^{\ast},\delta=\delta^{\ast},\lambda=\lambda^{\ast}$ or,
because of (\ref{404})-(\ref{406}),%
\begin{equation}
\kappa=\left \vert \delta \right \vert ^{3/2},\quad \lambda=\frac{2\delta \delta \,%
\acute{}%
\,%
\acute{}%
-3\delta \,%
\acute{}%
\,^{2}-4\delta^{2}}{4\left \vert \delta \right \vert ^{7/2}}. \tag{4.7}%
\label{407}%
\end{equation}
Thus we deduce

\begin{proposition}
Let $\Phi \subset%
\mathbb{R}
^{3}$ be a non-conoidal ruled $C^{3}$-surface, free of torsal rulings, which
is parametrized by (\ref{301}).\ Then $\Phi$ is congruent with its affine
normal image $\Phi^{\ast}$ if and only if the fundamental invariants of $\Phi$
are associated with (\ref{407}).
\end{proposition}

\noindent By virtue of (\ref{303}), we can easily confirm that the vectors
$\bar{s}$\thinspace$%
\acute{}%
$\thinspace$(u)$ and $\bar{r}^{\ast}$\thinspace$%
\acute{}%
$\thinspace$(u)$ are parallel or orthogonal if and only if $\delta \,%
\acute{}%
=1+\kappa \lambda=0$ or $\kappa=\lambda$ respectively. The surface $\Phi$ is an
Edlinger surface in the first case. In the second case the line of striction
of $\Phi$ is an asymptotic line. We formulate these results in the following

\begin{proposition}
Let $\Phi \subset%
\mathbb{R}
^{3}$ be a non-conoidal ruled $C^{3}$-surface, free of torsal rulings, and
$\Phi^{\ast}$ its affine normal image$.$ If the tangents of the line of
striction of $\Phi$ and the directrix $\Gamma$ of $\Phi^{\ast}$ in the
corresponding points are parallel (resp. orthogonal), then $\Phi$ is an
Edlinger surface (resp. the line of striction of $\Phi$ is an asymptotic line).
\end{proposition}

\noindent Let $\Phi^{\ast}$ be an Edlinger surface, i.e. $\delta^{\ast}\,%
\acute{}%
=1+\kappa^{\ast}\lambda^{\ast}=0$. On account of (\ref{404})-(\ref{406}) we
obtain $\kappa=const.$, $\delta=const.$ or%
\begin{equation}
\kappa=\frac{c_{0}}{u},\quad \delta=\frac{c_{1}}{u^{2}},\quad c_{0}%
,c_{1}=const.,\quad c_{0}c_{1}\neq0. \tag{4.8}\label{408}%
\end{equation}
The condition $\kappa=const.$ means that $\Phi$ is of constant slope. So we
can state

\begin{proposition}
Let $\Phi \subset%
\mathbb{R}
^{3}$ be a non-conoidal ruled $C^{3}$-surface, free of torsal rulings, which
is parametrized by (\ref{301}). If the affine normal image $\Phi^{\ast}$ of
$\Phi$ is an Edlinger surface, then $\Phi$ is a ruled surface of constant
slope and constant parameter of distribution or the relations (\ref{408}) are valid.
\end{proposition}

\noindent If both ruled surfaces $\Phi$ and $\Phi^{\ast}$ are Edlinger
surfaces, then their fundamental invariants are constant. Consequently we have

\begin{proposition}
Assume that the affine normal image of an Edlinger $C^{3}$-surface is an
Edlinger surface too. Then the osculating hyperboloids of each of them are congruent.
\end{proposition}

We conclude this work by studying the rulings-preserving mapping
$f:\Phi \longrightarrow \Phi^{\ast}$ defined by considering the parametrizations
(\ref{301}), (\ref{401}) and making the association $\bar{x}%
(u,v)\longrightarrow \bar{x}^{\ast}(u,v).$ The components $g_{ij}^{\ast}$ of
the first fundamental tensor of $\Phi^{\ast}$ in the coordinates $u^{1}:=u$,
$u^{2}:=v$ are the following%
\begin{equation}
g_{11}^{\ast}=\left \vert \delta \right \vert ^{-5}\left[  \frac{2\delta \delta \,%
\acute{}%
\,%
\acute{}%
-3\delta \,%
\acute{}%
\,^{2}-4\delta^{2}}{4}+\frac{\left(  2\delta \kappa \,%
\acute{}%
-3\delta \,%
\acute{}%
\, \kappa \right)  v}{2}\right]  ^{2}\!+\kappa^{2}v^{2}\left \vert
\delta \right \vert ^{-3}+\kappa^{2}\left \vert \delta \right \vert ^{-1},
\tag{4.9}\label{409}%
\end{equation}%
\begin{equation}
g_{12}^{\ast}=g_{21}^{\ast}=\varepsilon \kappa \left \vert \delta \right \vert
^{-4}\left[  \frac{2\delta \delta \,%
\acute{}%
\,%
\acute{}%
-3\delta \,%
\acute{}%
\,^{2}-4\delta^{2}}{4}+\frac{\left(  2\delta \kappa \,%
\acute{}%
-3\delta \,%
\acute{}%
\, \kappa \right)  v}{2}\right]  , \tag{4.10}\label{410}%
\end{equation}%
\begin{equation}
g_{22}^{\ast}=\kappa^{2}\left \vert \delta \right \vert ^{-3}. \tag{4.11}%
\label{411}%
\end{equation}
Let $f$ be an area-preserving mapping. Then $\det \left(  g_{ij}\right)
=\det(g_{ij}^{\ast})$, which, on account of (\ref{304}) and (\ref{409}%
)-(\ref{411}), is equivalent to%

\begin{equation}
\kappa^{4}\left \vert \delta \right \vert ^{-6}v^{2}+\kappa^{4}\left \vert
\delta \right \vert ^{-4}=v^{2}+\delta^{2}. \tag{4.12}\label{412}%
\end{equation}
Hence%
\begin{equation}
\kappa=\varepsilon_{0}\left \vert \delta \right \vert ^{3/2},\quad \varepsilon
_{0}=\pm1. \tag{4.13}\label{413}%
\end{equation}
Let the mapping $f$ be conformal. Then
\begin{equation}
\frac{g_{11}^{\ast}}{g_{11}}=\frac{g_{12}^{\ast}}{g_{12}}=\frac{g_{22}^{\ast}%
}{g_{22}}. \tag{4.14}\label{414}%
\end{equation}
Inserting (\ref{409})-(\ref{411}) in (\ref{414}) and taking into account
(\ref{304}), we find
\begin{equation}
\kappa=c\left \vert \delta \right \vert ^{3/2},\quad \lambda=\frac{2\delta \delta \,%
\acute{}%
\,%
\acute{}%
-3\delta \,%
\acute{}%
\,^{2}-4\delta^{2}}{4c\left \vert \delta \right \vert ^{7/2}},\quad
c=const.\neq0. \tag{4.15}\label{415}%
\end{equation}
Especially, the conformal mapping $f$ is an isometry (of Minding) if and only
if $c=\varepsilon_{0}.$ So we have the following result:

\begin{proposition}
The above mentioned mapping $f:\Phi \longrightarrow \Phi^{\ast}$ is an
area-preserving mapping (resp. conformal), if the fundamental invariants of
$\Phi$ are associated with (\ref{413}) (resp. (\ref{415})). In particular $f$
is an isometry if and only if $c=\varepsilon_{0}.$
\end{proposition}

\begin{remark}
By an isometry $f:\Phi \longrightarrow \Phi^{\ast}$ the fundamental invariants
($\kappa^{\ast},\delta^{\ast},\lambda^{\ast}$) of $\Phi^{\ast}$ are
($\kappa,\varepsilon_{0}\delta,\lambda$). This means that the ruled surfaces
$\Phi$~and~$\Phi^{\ast}$ are congruent or opposite congruent (in German
"gegensinnig kongruent") if $\varepsilon_{0}=+1$ or $\varepsilon_{0}=-1$
respectively \cite[p. 23]{Hoschek}.
\end{remark}

\noindent \textbf{Acknowledgement}

\noindent The authors would like to express their thanks to the referee for
his useful remarks.

\end{document}